\documentclass{article}
\usepackage{spconf,amsmath,amssymb,amsthm,epsfig,dsfont,color,subfigure,empheq,graphicx,subfig}
\usepackage{enumerate,url,algpseudocode,algorithm,wasysym,epstopdf}
\usepackage[table]{xcolor}

\newcommand \ba{\mathbf{a}}
\newcommand \bb{\mathbf{b}}

\newcommand \bone{\mathbf{1}}
\newcommand \bef{\mathbf{f}}

\newcommand \bv{\mathbf{v}}

\newcommand \by{\mathbf{y}}
\newcommand \bz{\mathbf{z}}
\newcommand \bp{\mathbf{p}}
\newcommand \bq{\mathbf{q}}

\newcommand \bC{\mathbf{C}}

\newcommand \bK{\mathbf{K}}

\newcommand \bQ{\mathbf{Q}}

\newcommand \bX{\mathbf{X}}
\newcommand \bY{\mathbf{Y}}

\newcommand \bR{\mathbf{R}}

\newcommand \tbq{\tilde{\mathbf{q}}}

\newcommand \mcH{\mathcal{H}}
\newcommand \mcK{\mathcal{K}}

\newcommand \mcN{\mathcal{N}}
\newcommand \mcQ{\mathcal{Q}}
\newcommand \mcT{\mathcal{T}}
\newcommand \mcZ{\mathcal{Z}}

\DeclareMathOperator{\trace}{Tr}

\newtheorem{lemma}{Lemma}

\begin{document}
\title{Kernel-Based Learning for Smart Inverter Control}


\twoauthors
 {Aditie Garg, Mana Jalali, Vassilis Kekatos\sthanks{This work is partially supported by grant NSF-CAREER-1751085.}}
	{Dept. of ECE, Virginia Tech, \\
		Blacksburg, VA 24061, USA}
  {Nikolaos Gatsis\sthanks{This work is partially supported by grant NSF-CCF-1421583.}}
	{Dept. of ECE, Un. of Texas at San Antonio, \\
		San Antonio, TX 78249, USA}
\maketitle

\begin{abstract}
Distribution grids are currently challenged by frequent voltage excursions induced by intermittent solar generation. Smart inverters have been advocated as a fast-responding means to regulate voltage and minimize ohmic losses. Since optimal inverter coordination may be computationally challenging and preset local control rules are subpar, the approach of customized control rules designed in a quasi-static fashion features as a golden middle. Departing from affine control rules, this work puts forth non-linear inverter control policies. Drawing analogies to multi-task learning, reactive control is posed as a kernel-based regression task. Leveraging a linearized grid model and given anticipated data scenarios, inverter rules are jointly designed at the feeder level to minimize a convex combination of voltage deviations and ohmic losses via a linearly-constrained quadratic program. Numerical tests using real-world data on a benchmark feeder demonstrate that nonlinear control rules driven also by a few non-local readings can attain near-optimal performance.
\end{abstract}


\begin{keywords}
Smart inverters, collaborative filtering, voltage regulation, power loss minimization.
\end{keywords}

\section{Introduction}\label{sec:intro}
While the stochasticity of distributed energy resources challenges grid operation, the smart inverters interfacing photovoltaic units can be engaged in reactive power compensation~\cite{Turitsyn11}. However, coordinating hundreds of inverters distributed over a feeder in real-time is non-trivial. 

Purely localized schemes utilizing local voltage or solar generation readings have been proposed for reactive power control by smart inverters~\cite{Turitsyn11}; yet their equilibria do not coincide with sought minimizers~\cite{FCL,LQD14,VKZG16}. In fact, there exist cases where local rules perform worse than the no-reactive support option~\cite{Jabr18}. Centralized grid dispatch schemes entail communicating electric loads and solar generation to the operator; solving an optimal power flow (OPF); and sending the computed setpoints back to inverters~\cite{FCL}, \cite{Robbins}, \cite{Rogers10}. Decentralized OPF solvers typically need several communication exchanges between neighboring inverters to converge~\cite{tse2014emiliano}, \cite{QiuyuPeng}, \cite{BaGa-TSG2017}. 

To reduce the computational and communication overhead, inverter control rules have been postulated as \emph{affine rules or policies}, evaluated at uncertain loads and generation~\cite{Jabr18}, \cite{Ayyagari17}, \cite{LinThomasBitar18}, \cite{Baker18}. Originally adopted under stochastic OPF setups in transmission systems~\cite{Bienstock}, an affine policy expresses a control variable as a linear function of given system inputs. The related weights are optimized periodically in a centralized fashion, yet the control rules are applied in real time. In distribution systems, controlling inverters via affine policies has been reported via chance-constrained~\cite{Ayyagari17}, robust \cite{Jabr18}, \cite{LinThomasBitar18}, and closed-loop formulations~\cite{Baker18}. Optimal policies however are not necessarily linear even for linear systems: If an inverter constraint becomes active, optimal power injections can become nonlinear functions of solar generation. We put forth data- and model-driven nonlinear policies using kernel-based learning. 


\emph{Notation:} lower- (upper-) case boldface letters denote column vectors (matrices). Calligraphic symbols are reserved for sets. Symbol $^{\top}$ stands for transposition, and $\mathbf{1}$ is the vector of all ones. The symbol $\|\mathbf{x}\|_2$ denotes the $\ell_2$-norm of $\mathbf{x}$, and $\|\cdot\|_F$ is the Frobenius matrix norm with $\|\bX\|_F^2=\trace(\bX^\top\bX)$.
 
\section{Problem Formulation}\label{sec:problem}
Before formulating the task of inverter control rule design, let us introduce the grid quantities needed. A single-phase radial grid with $N+1$ buses can be modeled by a tree graph whose nodes correspond to buses, and edges to distribution lines. The substation is indexed by $n=0$, and the remaining buses comprise the set $\mcN=\{1,\ldots,N\}$. Let $(v_n,p_n,q_n)$ denote the voltage magnitude and (re)active injections at bus $n$. 

The active power injected at bus $n$ can be decomposed as $p_n=p_n^g-p_n^c$, where $p_n^g$ is the solar generation and $p_n^c$ the inelastic load at the same bus. Reactive power injections can be similarly decomposed as $q_n=q_n^g-q_n^c$. For known $p_n^g$, the reactive injection of inverter $n$ is constrained by its apparent power limit $\bar{s}_n^g$ as
\begin{equation}\label{eq:pv}
|q_n^g|\leq \bar{q}_n^g:=\sqrt{(\bar{s}_n^g)^2 - (p_n^g)^2}.
\end{equation}
The $N$--length vectors $\bp=\bp^g-\bp^c$, $\bq=\bq^g-\bq^c$, $\bv$, and $\bar{\bq}^g$ collect the related nodal quantities and limits. We adopt the \emph{linearized distribution flow} (LDF) model, according to which voltages are approximately linear functions of power injections~\cite{BW3}, \cite{SGC15}
\begin{equation}\label{eq:LDF}
\bv \simeq \bR\bp + \bX\bq + v_0\mathbf{1}
\end{equation}
where $\mathbf{1}$ is the all-one vector; and matrices $(\bR,\bX)$ depend on the network topology and impedances~\cite{SGC15}.

Given solar generation and electric loads $(\bp^g,\bp^c,\bq^c)$, the task of \emph{reactive power compensation} by smart inverters aims at finding the setpoints $\bq^g$ to: \emph{i)} maintain voltage deviations within the range imposed by the ANSI C84.1 standard; and \emph{ii)} minimize the thermal (ohmic) losses on distribution lines. For the former objective, the operator may want to minimize the squared voltage deviations, which can be expressed as $V(\bq^g):=\|\bv-v_0\mathbf{1}\|_2^2\simeq \|\bR\bp+\bX\bq\|_2^2$ thanks to \eqref{eq:LDF}. Ohmic losses can be approximately expressed as $\bp^\top\bR\bp +\bq^\top\bR\bq$; see~\cite{Turitsyn11} for details. Since the control variable $\bq^g$ appears only in the second summand, ohmic losses simplify as
\begin{equation}\label{eq:losses}
L(\bq^g):=\bq^\top\bR\bq.
\end{equation}
The positive definiteness of $\bR$ guarantees that $L(\bq^g)$ is a positively-valued convex quadratic function. 

By definition, the objectives $V(\bq^g)$ and $L(\bq^g)$ are contradicting in general~\cite{Turitsyn11}. To handle this multiobjective optimization, the scalarization approach of \cite{Jabr18} can be adopted to pose the reactive power compensation task as 
\begin{equation}\label{eq:rpc1}
\min_{\bq^g\in \mcQ\nonumber}~\lambda V(\bq^g) + (1-\lambda) L(\bq^g)
\end{equation}
where the set $\mcQ:=[-\bar{\bq}^g,\bar{\bq}^g]$ captures the apparent power constraints in \eqref{eq:pv} for all $n\in\mcN$. By solving \eqref{eq:rpc1} for different values of $\lambda\in[0,1]$, the Pareto front for this control task can be recovered. Being a convex combination of two convex quadratic functions, the cost in \eqref{eq:rpc1} is apparently convex too. Upon completing the squares and ignoring inconsequential terms, problem \eqref{eq:rpc1} can be simplified as follows.

\begin{lemma}\label{le:transform}
Problem \eqref{eq:rpc1} can be equivalently expressed as
\begin{equation}\label{eq:rpc2}
\tbq^g:=\arg\min_{\bq^g\in \mcQ}~\|\bC\bq^g+\by\|_2^2
\end{equation}
where matrix $\bC:=[(1-\lambda)\bR+\lambda \bX^2]^{1/2}$; vector $\by:=\bC^{-1}[-(1-\lambda)\bR\bq^c +\lambda\bX\bR(\bp^g-\bp^c) - \lambda\bX^2\bq^c]$; and the operator $[\cdot]^{1/2}$ represents the unique square root of a symmetric positive definite matrix.
\end{lemma}

It is worth noticing that matrix $\bC$ depends only on the feeder, whereas vector $\by$ and the set $\mcQ$ depend on variable loads and solar generations collected in vector $\bz:=[(\bp^c)^\top~(\bq^c)^\top~(\bp^g)^\top]^\top$. 

An ideal control process entails the ensuing three steps:\\
\hspace*{1em}\textbf{\textit{S1)}} Each node communicates its load and solar data $(p_n^g,p_n^c,q_n^c)$ to the operator;\\
\hspace*{1em}\textbf{\textit{S2)}} the operator solves \eqref{eq:rpc2} given $\bz$; and\\
\hspace*{1em}\textbf{\textit{S3)}} the operator sends the optimal setpoints $\tbq^g$ to inverters.

For varying $\bz$, this control process should be repeated on a per-minute basis or more frequently. Observe that step \textit{S1)} requires $N$ inverter-utility communication links, and \textit{S3)} another $N$ utility-inverter links. Running this operation for multiple feeders hosting hundreds of buses each, constitutes a computation- and communication-wise formidable task. To reduce this cyber overhead, the operator may decide to issue setpoints less frequently; but then setpoints may become obsolete and suboptimal. 

To bypass this limitation, we suggest designing \emph{control policies} according to which the reactive power injection from inverter $n$ is a function of grid data as
\begin{equation}\label{eq:qfun}
q_n^g(\bz_n)=f_n(\bz_n)+b_n
\end{equation}
where $f_n$ is an inverter-customized function; its argument $\bz_n\in\mcZ_n\subseteq \mathbb{R}^{M_n}$ is a subvector of grid data $\bz$; and $b_n$ is an intercept. For a purely local rule, the control input can be selected as 
\begin{equation*}
\bz_n:=[p_n^g~~\bar{q}_{n}^g~~p_n^c~~q_n^c]^\top. 
\end{equation*}
If communication resources are abundant, one can set $\bz_n=\bz$. Otherwise, hybrid scenarios could be obviously envisioned.

Designing inverter rules via control policies has been advocated in~\cite{Jabr18}, \cite{Ayyagari17}, \cite{LinThomasBitar18}, \cite{Baker18}. Yet the control policies were confined to linear $f_n$'s. Upon reviewing kernel-based learning, the next section provisions \emph{non-linear} mappings $f_n$.

\section{Learning Inverter Control Rules}\label{sec:background}
Kernels have served as the foundation for extending machine learning tools to nonlinear mappings. Given pairs $\{(\bz_t,y_t)\}_{t=1}^T$ of features $\bz_t$ belonging to a space $\mcZ$ and target values $y_t\in \mathbb{R}$, kernel-based learning aims at finding a mapping $f:\mcZ\rightarrow\mathbb{R}$. The mapping $f$ is constrained to lie on the linear function space~\cite{Evgeniou00}
\begin{equation}\label{eq:family}
\mcH_\mcK:=\left\{f(\bz)=\sum_{t=1}^{\infty} K(\bz,\bz_t) a_t,~a_t\in\mathbb{R}\right\}
\end{equation}
defined by a \emph{kernel function} $K:\mathcal{Z}\times \mathcal{Z}\rightarrow \mathbb{R}$ and coefficients $a_t$. When $K(\cdot,\cdot)$ is a symmetric positive definite function, the function space $\mathcal{H}_{\mathcal{K}}$ becomes a reproducing kernel Hilbert space (RKHS) whose members have a finite norm~\cite{Evgeniou00}, \cite{BaGia13}
\[\|f\|_{\mathcal{K}}^2:= \sum_{t=1}^{\infty} \sum_{t'=1}^{\infty} K(\bz_t,\bz_{t'}) a_t a_{t'}.\]

Learning $f$ from data $\{(\bz_t,y_t)\}_{t=1}^T$ can be formulated as the functional minimization task~\cite{BaGia13}
\begin{equation}\label{eq:fnreg}
\hat{f}:=\arg\min_{f} \frac{1}{T}\sum_{t=1}^T \left[y_t-f(\bz_t)\right]^2 + \mu \|f\|_{\mathcal{K}}^2.
\end{equation}
The first summand in \eqref{eq:fnreg} is a data-fitting term. The second one ensures that $\|\hat{f}\|_{\mathcal{K}}$ is finite, and so $\hat{f}\in\mcH_{\mcK}$. More complex functions have higher $\|f\|_{\mathcal{K}}$, yield a better fit to training data, but can perform poorly on unseen data. The parameter $\mu>0$ balances fitting over generalization and can be tuned via cross-validation. 

The celebrated Representer's Theorem asserts that the minimizer of \eqref{eq:fnreg} takes the form $\hat{f}(\bz)=\sum_{t=1}^T K(\bz,\bz_t)\hat{a}_t$, that is $\hat{f}$ is described only by $T$ rather than infinitely many coefficients $a_t$'s. Then, the objective of \eqref{eq:fnreg} can be expressed in terms of the unknown $a_t$'s, and the functional minimization is converted to a quadratic optimization.

Returning to the task of designing inverter control rules, the idea here is to leverage kernel-based learning and postulate that the mapping $f_n$ for inverter $n$ in \eqref{eq:qfun} lies in the RKHS
\begin{equation}\label{eq:rpcfamily}
\mcH_{\mcK_n}:=\left\{f_n(\bz_n)=\sum_{t=1}^{\infty} K_n(\bz_{n},\bz_{n,t}) a_{n,t},~a_{n,t}\in\mathbb{R}\right\}
\end{equation}
defined by the kernel function $K_n:\mcZ_n\times \mcZ_n\rightarrow \mathbb{R}$ with control inputs $\bz_n\in\mcZ_n$. Linear policies can be captured by selecting the linear kernel $K_n(\bz_{n,t},\bz_{n,t'})=\bz_{n,t}^\top\bz_{n,t'}$. Nonlinear policies can be designed by selecting for example a polynomial kernel $K_n(\bz_{n,t},\bz_{n,t'})=\left(\bz_{n,t}^\top\bz_{n,t'}+\gamma\right)^\beta$, or a Gaussian kernel $K_n(\bz_{n,t},\bz_{n,t'})=\exp\left(-\|\bz_{n,t}-\bz_{n,t'}\|_2^2/\gamma\right)$ with design parameters $\beta$ and $\gamma>0$. 

We propose designing control rules using $T$ scenario data $\left\{\bz_{n,t}\right\}_{n\in\mcN}$ over $t\in\mcT:=\{1,\ldots,T\}$, by solving the functional minimization
\begin{subequations}\label{eq:funmin}
\begin{align}
\min~&~\frac{1}{T}\sum_{t=1}^T C_\lambda\left[\left\{q_n^g(\bz_{n,t})\right\};\by_t\right] +\mu \sum_{n=1}^N\|f_n\|_{\mcK_n}^2\label{eq:rpcfnreg:cost}\\
\mathrm{over}~&~ q_n^g(\bz_{n,t})\in \mcH_{\mcK_n}~n\in\mcN,\quad \bb\in\mathbb{R}^N\\
\mathrm{s.to}~&~\eqref{eq:qfun},~|q_n^g(\bz_{n,t})|\leq \bar{q}_{n,t}^g,\quad\quad n\in\mcN,~t\in\mcT.\label{eq:rpcfnreg:con}
\end{align}
\end{subequations}
The least-square fit of \eqref{eq:fnreg} has been replaced by the cost of Lemma~\ref{le:transform} averaged over the training scenarios
$C_\lambda\left[\bq^g(\bz);\by\right]:=\|\bC\bq^g(\bz)+\by\|_2^2$.
Similar to collaborative filtering~\cite{KZG14,BaGia13}, the minimization in \eqref{eq:funmin} intends to learn $N$ rather than one function $\{q_n^g\}_{n\in\mcN}$. An alternative approach for learning control rules has appeared in \cite{Dobbe18}: Inverter rules are trained upon fitting directly control inputs to optimal inverter decisions. Different from our approach where inverter functions are naturally coupled through the underlying physical system, the scheme in \cite{Dobbe18} treats the electric grid and the OPF solver as a black box and trains each inverter function independently.

The control process is organized into four steps:

\textbf{\textit{T1) Data collection.}} On a 30-min basis, the operator collects smart meter readings (active and reactive loads, solar generation) from all buses. These readings and possibly historical data can be used as training scenarios $\{\bz_{n,t}\}$ for all $n\in\mcN$ and $t\in\mcT$.

\textbf{\textit{T2) Control rule design.}} On a 30-min basis, the operator finds simultaneously the control rules for all inverters by solving~\eqref{eq:funmin}. Fortunately, the functional minimization in \eqref{eq:funmin} can be converted to a vector optimization problem as elaborated in the next lemma. Before doing so, let us define the $T\times T$ kernel matrix $\bK_n$ for inverter $n$ with $(t,t')$ entry being equal to the kernel function evaluation $K_n(\bz_{n,t},\bz_{n,t'})$.

\begin{lemma}\label{le:LS}
The functional minimization in \eqref{eq:funmin} can be equivalently expressed as 
\begin{subequations}\label{eq:krpc}
\begin{align}
\min~&~\frac{1}{T}\|\bC\bQ+\bY\|_F^2 + \mu \sum_{n=1}^N \ba_n^\top\bK_n\ba_n \label{eq:krpc:cost}\\
\mathrm{over}~&~\bQ\in\mathbb{R}^{N\times T}, \{\ba_n\in\mathbb{R}^T\}_{n=1}^N, \bb\in\mathbb{R}^N \label{eq:krpc:vars}\\
\mathrm{s.to}~&~\bQ^\top=\left[ \bK_1\ba_1+b_1\bone~~\cdots~~ \bK_N \ba_N+b_N\bone\right]\label{eq:krpc:con1}\\
~&~-\bar{\bq}_n^g\leq \bK_n\ba_n+b_n\bone \leq \bar{\bq}_n^g,~ \forall n\label{eq:krpc:con2}
\end{align}
\end{subequations}
where $\bY:=[\by_1~\cdots~\by_T]$ and the entries of vector $\bar{\bq}_n^g:=[\bar{q}_{n,1}^g~\cdots~\bar{q}_{n,T}^g]^\top$ have been defined in \eqref{eq:pv}. 
\end{lemma}

\begin{proof}
Fortunately, the Representer's Theorem can be applied successively over $n$ in \eqref{eq:funmin}. It can thus ensure that the $f_n$ minimizing \eqref{eq:funmin} has the form
\begin{equation}\label{eq:funRT}
f_n(\bz_n)=\sum_{t=1}^T K_n(\bz_n,\bz_{n,t}) a_{n,t}
\end{equation}
for all $n$. Evaluating the inverter policy $f_n$ of \eqref{eq:funRT} over the test data $\{\bz_{n,t}\}_{t=1}^T$ yields for all $n$
\begin{equation*}
\bef_n = \bK_n \ba_n
\end{equation*}
where $\ba_n:=[{a}_{n,1}~\cdots~{a}_{n,T}]^\top$. The reactive power injections for inverter $n$ over all scenarios in $\mcT$ can then be expressed as
\begin{equation}\label{eq:bqng}
\bq_n^q=\bK_n\ba_n+b_n\bone.
\end{equation}
Then, the linear inequalities in \eqref{eq:krpc:con2} capture the apparent power constraints evaluated over the tested scenarios. 

Moreover, from the reproducing properties of $\mcH_{\mcK_n}$'s, the RKHS norms in the second summand of \eqref{eq:funmin} can be written as
\begin{equation*}
\|f_n\|_{\mcK_n}^2=\ba_n^\top\bK_n\ba_n,\quad \forall n.
\end{equation*}

Consider finally the first summand in \eqref{eq:krpc:cost}. Based on \eqref{eq:bqng}, the $t$-th column of $\bQ$ denoted by $\bq_t^g$ contains the reactive injections from all inverters at scenario $t$. Then, the data-fitting term of \eqref{eq:funmin} can be written as
\[\sum_{t=1}^T\|\bC\bq_t^g+\by_t\|_2^2=\|\bC\bQ+\bY\|_F^2\]
since the squared Frobenius norm of a matrix equals the sum of the squared $\ell_2$-norms of its columns.
\end{proof}

Lemma~\ref{le:LS} poses the task of finding optimal control policies in \eqref{eq:funmin} as a linearly-constrained quadratic program over $\{\ba_n,b_n\}$ for all $n\in\mcN$. Once the latter parameters have been found via \eqref{eq:krpc}, the policy mappings $\{f_n\}$ can be evaluated for any other input $\bz_n$ using \eqref{eq:funRT}.

\textbf{\textit{T3) Downloading control rules.}} Upon solving \eqref{eq:krpc}, the rule for inverter $n$ is fully described by $(\ba_n,b_n)$ and scenarios $\{\bz_{n,t}\}_{t\in\mcT}$. If $\bz_{n,t}\in\mathbb{R}^{M_n}$, the operator needs to send $(M_n+1)T+1$ data to inverter $n$. However, this step occurs once every 30~min.

\textbf{\textit{T4) Control rule implementation.}} In near real-time (say every 30 sec) and for the next 30 min, each inverter $n$ implements its control rule of \eqref{eq:qfun} by evaluating the kernel function $K_n(\bz_n,\bz_{n,t})$ for the feeder conditions $\bz_n$ currently experienced. If the control input $\bz_n$ is purely local, no communication is needed. Otherwise, remote inputs have to be communicated from their sources to inverter $n$. 

Although the constraints in \eqref{eq:rpcfnreg:con} are enforced for the test data, the policies obtained via \eqref{eq:funmin} may not satisfy the apparent power limits for $\bz_{n,t}$'s with $t\notin \mcT$. This limitation of kernel-based learning appears also in scenario-based and chance-constrained designs~\cite{Ayyagari17}. Of course, once the extrapolated control policy has been found from \eqref{eq:qfun}, its value can be heuristically projected within $[-\bar{q}_{n,t'}^g,+\bar{q}_{n,t'}^g]$ as
\begin{equation*}\label{eq:project}
\left[q_n^g(\bz_{n,t})\right]_{\bar{q}_{n,t}^g}:=\max\left\{\min\left\{q_n^g(\bz_{n,t}),\bar{q}_{n,t}^g\right\},-\bar{q}_{n,t}^g\right\}.
\end{equation*}

In the standard machine learning setup, one wants to fit a function to capture the dependency between features and targets. Ideally, the designed function should behave well even for feature-target pairs not seen during training or fitting process. In direct analogy, the inverter control policies are posed as a joint function fitting task based on scenario data. The grid quantities feeding each controller serve as feature data, and the reactive injections as target values. Once the functions have been designed, they can be applied to unseen data.

\section{Numerical Tests}\label{sec:tests}
\begin{figure}[t]
\centering
\includegraphics[scale=0.28]{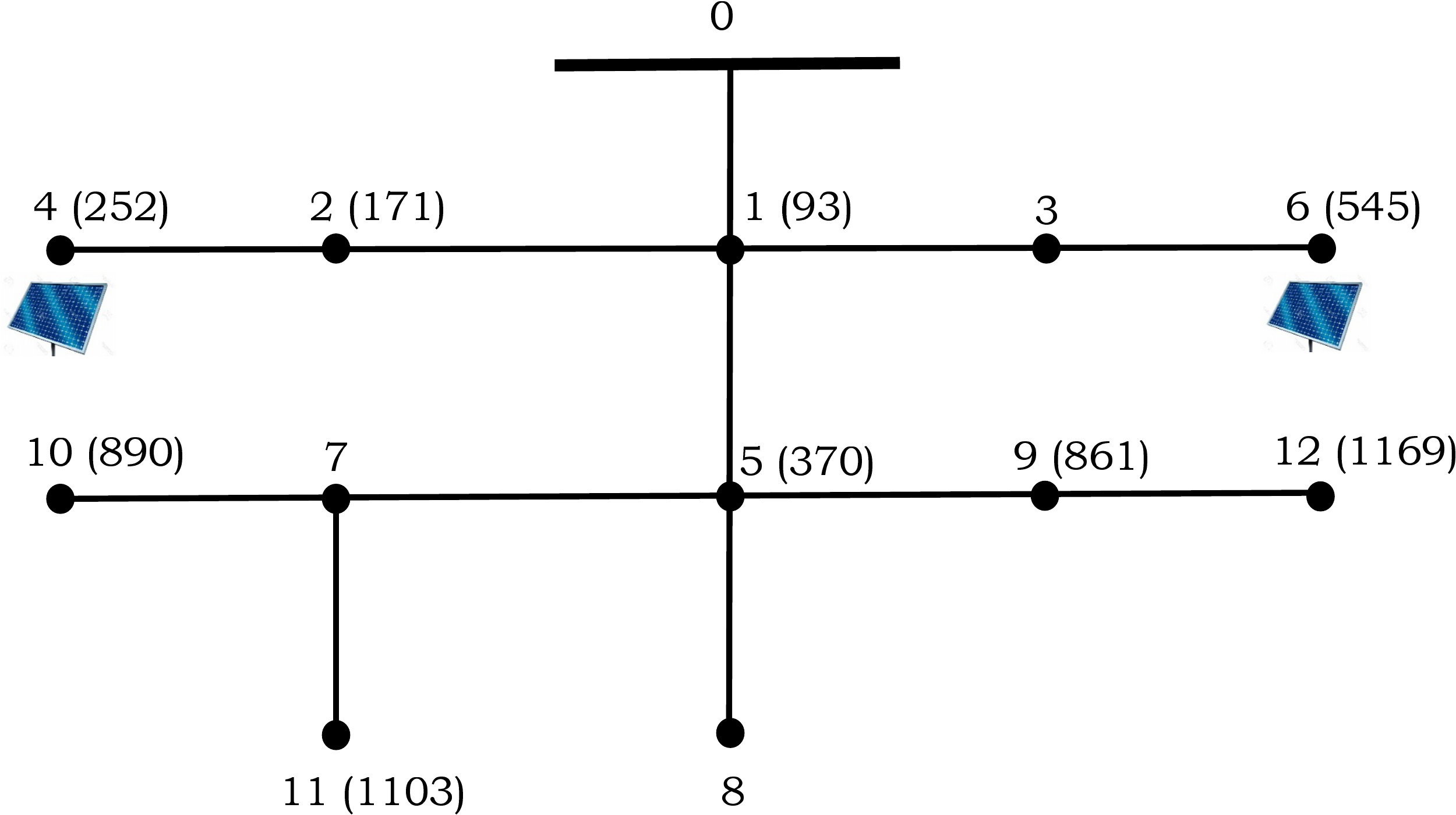}
\caption{Matching between buses and Pecan Str houses.}
\label{fig:13bus}
\end{figure}

Our control rules were tested on the IEEE 13-bus feeder, converted to a single-phase grid~\cite{GLTL12}. Minute-sampled loads and solar generation data were extracted from the Pecan Str dataset at \url{https://dataport.cloud/} for October 1, 2013. Figure~\ref{fig:13bus} matches the Pecan Str house indexes to nonzero-injection buses. Reactive loads were randomly drawn to yield power factors uniformly distributed in $[0.90,0.95]$ lagging. Each load timeseries was scaled so that its monthly peak matched $50\%$ of the benchmark load. Solar data were scaled using the previous numbers. 


\begin{figure}[t]
\centering
\includegraphics[scale=0.30]{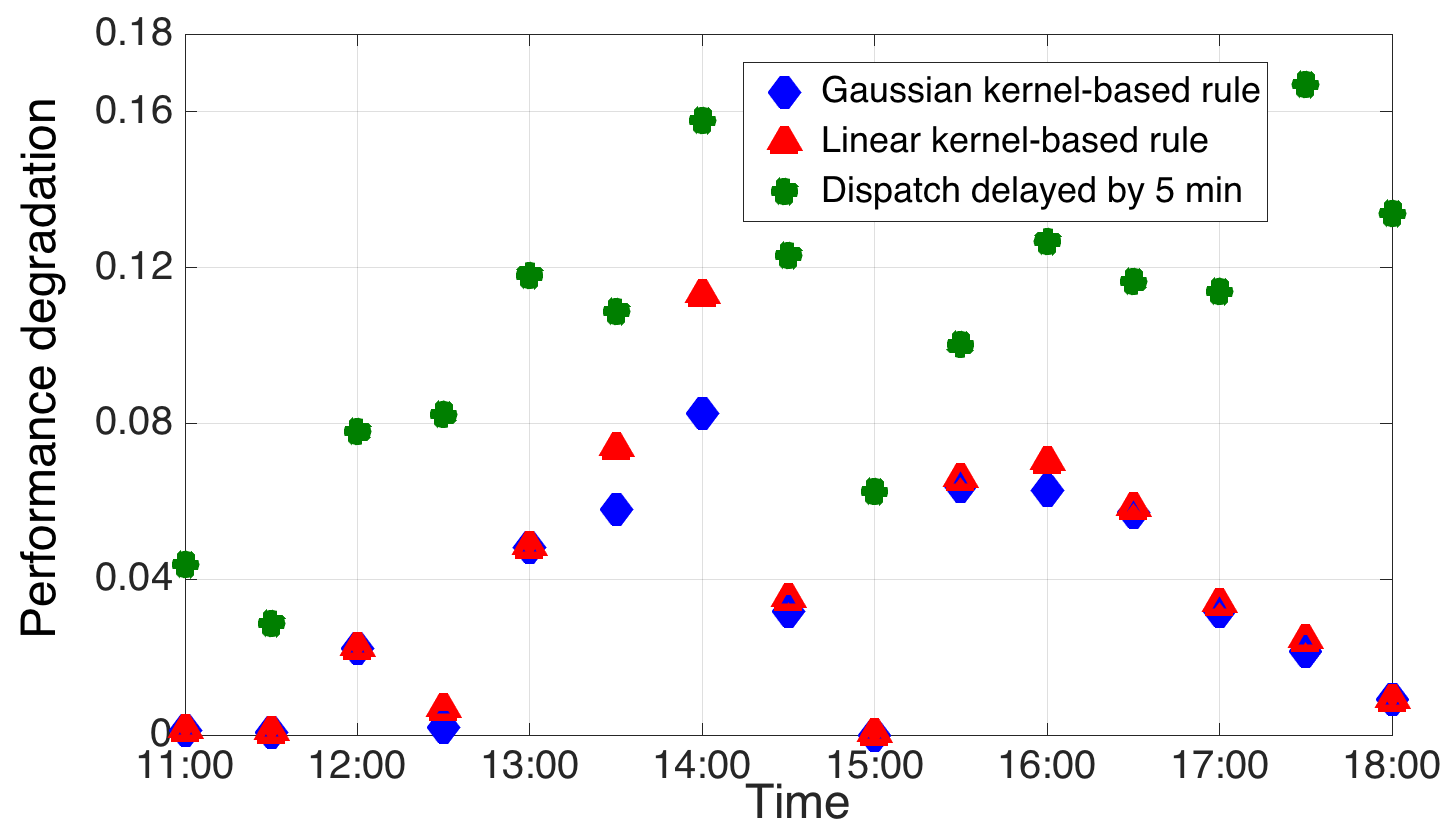}
\caption{Performance degradation to the optimal dispatch averaged over 30-min intervals for $\lambda=0.5$.}
\label{fig:cost}
\end{figure}

We tested the optimal reactive power setpoints obtained via \eqref{eq:rpc2} on a per-minute basis, and the kernel-based rules of \eqref{eq:krpc} for the linear and Gaussian kernels. The rules were trained using $T=30$ minutes or scenarios, while $\mu$ and $\gamma$ were set using 5-fold cross-validation. The controllers' input comprised local data along with the active flows on lines $(1,2)$, $(1,3)$, and $(1,5)$. Problem \eqref{eq:rpc2} was solved using the MATLAB-based toolbox YALMIP along with the SDPT3 solver~\cite{yalmip}, \cite{sdpt3}. The kernel-based control rules were solved using the OSQP solver~\cite{osqp}, and applied over the next 30 minutes. For each scheme, the cost of \eqref{eq:Clambda} was evaluated for $\lambda=0.5$ and averaged over the $30$-min period between hours 11:00--18:00. Figure~\ref{fig:cost} shows the difference between the cost obtained by the rules and the optimal cost of~\eqref{eq:rpc2}. At times of low solar irradiance, the rules coincide with the optimal dispatch. During higher solar generation, the suboptimality of rules increases as expected, yet the Gaussian kernel-based rule outperforms the linear rule in general. To test the effect of obsolete optimal setpoints, the setpoints found via \eqref{eq:rpc2} for minute $t$ were applied to the system at time $t+5$. As with control rules, the setpoint $q_{n,t}^g$ was projected to comply with \eqref{eq:pv} for the current $p_{n,t+5}^g$. The performance degradation is significant over the one provided by rules. 
	
\section{Conclusions}\label{sec:conclusions}
Non-linear control policies for determining inverter reactive injections have been designed using the powerful tool of kernel-based learning. The policies are centrally designed on a 30-min basis, although they are run in real-time using local and/or remote grid data. The design is extremely flexible: its communication needs depend on local and remote inputs, while the computationally demanding task of \eqref{eq:krpc} is run at the utility on a 30-min basis. Numerical tests using real-world data on a benchmark feeder validate that the suggested control rules can achieve the desirable trade-off between feeder performance and cyber overhead.  

\vfill\pagebreak

\bibliographystyle{IEEEtran}
\bibliography{myabrv,globalsip}
\end{document}